\documentclass[12pt, a4paper]{article}
\usepackage{amssymb}
\usepackage{amsmath}
\usepackage[margin= 2.0cm]{geometry}
\usepackage{titlesec}
\usepackage{color}
\usepackage{graphicx}
\usepackage{subfig}
\usepackage{multirow}
\usepackage{array}

\newcommand{\comm}[1]{}

\titleformat{\section}{\Large\bfseries}{\thesection.}{.5em}{}
\titlespacing*{\section}{0pt}{10pt}{5pt}
\titleformat{\subsection}{\large \bfseries}{\thesubsection.}{.5em}{}
\titlespacing*{\subsection} {0pt}{10pt}{5pt}
\titleformat{\subsubsection}{\normalfont\bfseries}{\thesubsubsection.}{.5em}{}
\titlespacing*{\subsubsection} {0pt}{10pt}{5pt}
\titleformat{\paragraph}[runin]{\normalfont\it\bfseries}{\theparagraph.}{.5em}{}[.]
\titlespacing{\subsubsubsection} {0pt}{6pt}{*1}

\newcounter{ew}

\newcommand{\Var}{\operatorname{Var}}

\DeclareMathOperator*{\argmax}{arg\,max}
\input{tcilatex}
\begin{document}

\begin{center}
\textbf{On a representation of fractional Brownian motion and the limit
distributions of statistics arising in cusp statistical models}
\end{center}

\bigskip

\begin{center}
Nino Kordzakhia ({\emph{Macquarie University, Sydney}})\footnote{%
Macquarie University, Sydney, Australia; e-mail:Nino.Kordzakhia@mq.edu.au}

Yury Kutoyants ({\emph{Le Mans University, France}})\footnote{%
University of Maine, Le Mans, France, and National Research University
\textquotedblleft MPEI\textquotedblright , Moscow, Russia.
\par
Research of Y.Kutoyants was partially supported by the grant of RSF number
14-49-00079.}

Alex Novikov ({\emph{University of Technology Sydney}})\footnote{%
Steklov Mathematical Institute of RAS, Moscow, and University of Technology
Sydney. Present address: PO Box 123, Broadway, Department of Mathematical
Sciences, University of Technology, Sydney, NSW 2007, Australia;
e-mail:Alex.Novikov@uts.edu.au
\par
Research of A.Novikov was partially supported by ARC Discovery grant
DP150102758.}

Lin-Yee Hin ({\emph{University of Technology Sydney}})\footnote{%
University of Technology Sydney. Present address: PO Box 123, Broadway,
Department of Mathematical Sciences, University of Technology, Sydney, NSW
2007, Australia; e-mail: LinYee.Hin@uts.edu.au
\par
Research of L. Hin was supported by ARC Discovery grant DP150102758.}

\bigskip
\end{center}

\textbf{Abstract. }We discuss some extensions of results from the recent
paper by Chernoyarov et al. (Ann. Inst. Stat. Math., October 2016)
concerning limit distributions of Bayesian and maximum likelihood estimators
in "signal plus white noise" model with irregular cusp-type signals. Using a
new representation of fractional Brownian motion (fBm) in terms of cusp
functions we show that, as the noise intensity tends to zero, the limit
distributions can be expressed in terms of fBm for the full range of
asymmetric cusp-type signals correspondingly with the Hurst parameter $H$, $%
0<H<1$. The \ simulation results for the densities and variances of the
limit distributions of Bayesian and maximum likelihood estimators are also
provided.

\bigskip

\textbf{1. Introduction and main results. }The monograph of Ibragimov and
Khasminskii \cite{1} contains a powerful technique for studying asymptotic
properties of Bayesian estimators (BE) $\tilde{\theta}_{n}$ and maximum
likelihood estimators (MLE) $\hat{\theta}_{n}$ of a parameter $\theta ~\ $%
based on independent identically distributed (i.i.d.)\textbf{\ }observations 
$X^{n}=\left( X_{1},\ldots ,X_{n}\right) \ $with the marginal density
function $f\left( x,\theta \right) $. In particular, for irregular
statistical models they showed (see \cite{1}, Chapter 6, Theorem 6.2 and
Theorem 6.4 ) that the limit distributions of $\ \tilde{\theta}_{n}$ and$~%
\hat{\theta}_{n}$,\ as$\ \ n\rightarrow \infty $,$~$can be represented using
Poisson or Gaussian processes which in their turn are defined in terms of
the singularity points of a density function. The particular case of
cusp-type densities%
\begin{equation}
f\left( x,\theta \right) ={h(x,}{\theta )}\exp \left\{ -g(x,\theta
)|x-\theta |^{\alpha }\right\} ,~\theta \in \Theta =(\theta _{1},\theta
_{2}),~x\in R=(-\infty ,\infty ),  \label{example 1}
\end{equation}%
where $\alpha >0,~h~$and $g~$are smooth functions, was discussed in the
original paper \cite{IKh1974}, see also Chapter 6 in \cite{1}.

The question about efficiency of MLE\ in irregular i.i.d. statistical
experiments, in particular, with $\alpha >\frac{1}{2}$ in (\ref{example 1})
was raised by H. Daniels \cite{Daniels} who showed that the MLE is
asymptotically efficient and normal in this case. Subsequently, P. Rao\ \cite%
{Prakasa} showed that the limit distribution of $\hat{\theta}_{n}$ for $%
\alpha \in (0,1/2)~$can be expressed in terms of fractional Brownian motion
(fBm) with the Hurst parameter $H=\alpha +1/2\in (1/2,1)$~although the
question about its efficiency had not been addressed in \cite{Prakasa}.

Recall that continuous Gaussian process $W^{H}=\{W_{u}^{H},u\in R\}~$with$%
\,W_{0}^{H}=0,\ E(W_{u}^{H})=0$ is said to be a standard fBm with the Hurst
parameter $\,H\in (0,1]\ $if%
\begin{equation}
\,E|W_{u}^{H}-W_{s}^{H}|^{2}=|u-s|^{2H},\,\ u\in R,~s\in R.  \label{fBm}
\end{equation}%
A standard two-sided Brownian motion $W=(W_{u}^{1/2},~u\in R)$ is a
particular case of this definition.

Further we use the following notations:%
\begin{equation*}
L(\theta ,X^{n}):=\prod\limits_{i=1}^{n}f\left( X_{i},\theta \right)
\end{equation*}%
for the likelihood function;%
\begin{equation*}
\tilde{\theta}_{n}=\frac{\int_{_{_{-\infty }}}^{^{^{\infty
}}}u\,L(u,X^{n})q(u)du}{\int_{_{_{-\infty }}}^{^{^{\infty
}}}\,L(u,X^{n})q(u)du}
\end{equation*}%
for the BE with respect to quadratic loss function and the prior
distribution $q(\theta )$.

Set 
\begin{equation*}
Z_{u}^{(H)}:=\exp \{W_{u}^{H}-\frac{|u|^{2H}}{2}\},~u\in R.
\end{equation*}

Under some mild assumptions on $q(\theta )$ (including the case $q(\theta
)=1~$i.e. Pitman-type estimators, see \cite{Pit})$~$the theory developed in 
\cite{1} implies the following result for the i.i.d. cusp model (\ref%
{example 1}) with $\alpha \in (0,1/2),~H=\alpha +1/2$: 
\begin{equation}
(\tilde{\theta}_{n}-\theta )n^{^{\frac{1}{2H}}}/C_{H}\overset{d}{%
\longrightarrow }\zeta _{H}:=\frac{\int_{_{_{-\infty }}}^{^{^{\infty
}}}uZ_{u}^{(H)}du}{\int_{_{_{-\infty }}}^{^{^{\infty }}}Z_{u}^{(H)}du},
\label{limit1}
\end{equation}%
where $C_{H}$ is a known constant, the convergence $\overset{d}{%
\longrightarrow }$ is understood in distribution.

Furthermore, for MLE $\hat{\theta}_{n}$ it was shown in \cite{Prakasa} that%
\footnote{%
The uniqueness of $\xi _{H}$ with probability 1 is shown in \cite{Pflug}.} 
\begin{equation}
(\hat{\theta}_{n}-\theta )n^{^{\frac{1}{2H}}}/C_{H}\overset{d}{%
\longrightarrow }\xi _{H}:=\arg \max_{u\in R}(Z_{u}^{(H)}).
\label{MLE limit}
\end{equation}%
Hence, both BE $\tilde{\theta}_{n}$ and MLE $\hat{\theta}_{n}$~have the same
rate of convergence $n^{^{-\frac{1}{2H}}},~H\in (1/2,1)$ for the i.i.d. cusp
model (\ref{example 1}). Note that some general properties of $\zeta _{H},\
H\in (0,1),$ have been studied in \cite{NK}, \cite{NKL} where, in
particular, the positive finite constant $\lambda _{H}$ \ found is such that
for all $\lambda <\lambda _{H}$%
\begin{equation*}
E\exp \{\lambda |\zeta _{H}|^{2H}\}<\infty
\end{equation*}%
implying finiteness of the moments of $|\zeta _{H}|.$

In a similar context other continuous and discrete time models with fBm $%
W^{H}$ arising in the limits have been discussed in the monograph by
Kutoyants \cite{Kut book}, Dachian \cite{Datch 2003}, Gushchin and K\"{u}%
chler \cite{Gu}, D\"{o}ring \cite{Doring} and the references therein.~The
only paper, where the limits similar to (\ref{limit1}) and (\ref{MLE limit})
appear with $H\in (0,\frac{1}{2})$, is \cite{F10}, where an observed
diffusion process had the drift of the form $a|X_{t}-\theta |^{\alpha }$.
Note that in \cite{F10} it is assumed that the observed diffusion process is
a weak solution of the stochastic differential equation; however, for
defining of the likelihood ratio process the existence of a strong solution
is required and this fact had not been addressed.

In engineering and statistical literature there is a great interest to the
"signal plus white noise" type models, where observations $X^{T}=\left(
X_{t},0\leq t\leq T\right) \ $have the following dynamics 
\begin{equation}
dX_{t}=S(t,\theta )dt+\varepsilon dw_{t},\quad X_{0}=0,\quad 0\leq t\leq
T,~\varepsilon >0,~0<\theta _{1}<\theta <\theta _{2}<T.  \label{swn}
\end{equation}%
Here we assume that $w=\{w_{t},\ 0\leq t\leq T\}$ is a standard one-sided
Brownian motion, $S(t,\theta )$ is a "deterministic signal" which depends on
a parameter $\theta $ \ to be estimated, the \textquotedblleft finite
energy\textquotedblright\ condition 
\begin{equation}
\int_{0}^{T}S^{2}(u,\theta )du<\infty   \label{SS}
\end{equation}%
holds and $T$ is fixed. The important scenario \textit{\textquotedblleft
large signal-to-noise ratio\textquotedblright } corresponds to $\varepsilon
\rightarrow 0$.

The case of the observations $X^{nT}=\left( X_{t},~0\leq t\leq nT\right) $
of \ $T$-periodic signal $S(t,\theta )$ with $\varepsilon =1,~n=1,...,$ in (%
\ref{swn}) can be reduced to this model with $\varepsilon =\frac{1}{\sqrt{n}}%
,~X^{T}=\left( X_{t},\ 0\leq t\leq T\right) $ if we put 
\begin{equation*}
X_{t}=\frac{1}{n}\sum_{j=1}^{n}\left( X_{T\left( j-1\right) +t}-X_{T\left(
j-1\right) }\right) ,\ 0\leq t\leq T.
\end{equation*}

The model (\ref{swn}) is very different from (\ref{example 1}), however,
Chernoyarov et al. \cite{Chernoyar} showed that for cusp-type signals of the
form%
\begin{equation}
S(t,\theta )=a|t-\theta |^{\alpha },\ \alpha \in (0,\frac{1}{2}),
\label{CherDK}
\end{equation}%
the limit distributions of BE $\tilde{\theta}_{\varepsilon }$ and MLE $\hat{%
\theta}_{n}$ are also expressed in terms of fBm with ~$H=\alpha +1/2\in
(1/2,1)$.

One of the main results of our paper is the extension of results of \cite%
{Chernoyar} for model (\ref{swn}) to the case $H\in (0,1)$ under more
general assumptions, namely, when \textquotedblleft
signal\textquotedblright\ might be the asymmetric cusp function of the
following form\footnote{%
In \cite{Chernoyar} only symmetric cusp i.e. $a=b$ in (\ref{cusp1}) was
discussed.
\par
{}} 
\begin{equation}
S(t,\theta )=q_{\alpha }(t-\theta )+h(t,\theta ),  \label{signal}
\end{equation}%
\begin{equation}
q_{\alpha }(x):=(aI\{x\geq 0\}+bI\{x<0\})|x|^{\alpha },\ \alpha \in (-\frac{1%
}{2},\frac{1}{2}),\ a\geq 0,\ b\geq 0,\ ab\neq 0,  \label{cusp1}
\end{equation}%
where $h(t,\theta )\in C^{1,1}$ \ is a continuously differentiable with
respect both $t$ and $\theta ,$ $I\{.\}$ is the indicator function. Such
signals are unbounded for $\alpha \in (-\frac{1}{2},0);$\ however,$\ $%
condition (\ref{SS}) still holds. In the case of \ \textquotedblleft
discontinuous signals\textquotedblright\ from (\ref{swn}) with%
\begin{equation*}
S(t,\theta )=aI\{t-\theta >0\},~0\leq t\leq T,
\end{equation*}%
as it was shown \ in Section 7 of \ \cite{1}, the asymptotic results of the
type (\ref{limit1}) and (\ref{MLE limit}) hold with $H=1/2.$

Further we use the following notation with $q_{\alpha }$ from (\ref{cusp1})
and set 
\begin{equation}
\Gamma _{\alpha }^{2}:=\int_{-\infty }^{\infty }(q_{\alpha }(y-1)-q_{\alpha
}(y))^{2}dy,~\alpha \in (-\frac{1}{2},\frac{1}{2}).  \label{GGa}
\end{equation}%
Here and below we consider all stochastic integrals in the Ito's sense.

The following new representation for fBm will play an important role in
finding the limit distributions for the cusp model (\ref{swn}).

\textbf{Theorem 1}.\textit{\ Let }$W=(W_{y},\ y\in R)$\textit{\ be a
standard two-sided Bm, }$q_{a}(x)$\textit{\ be the cusp function from (\ref%
{cusp1}).~Then the process }%
\begin{equation}
Y_{u}^{H}:=\Gamma _{\alpha }^{-1}\int_{-\infty }^{\infty
}(q_{a}(y-u)-q_{a}(y))dW_{y},\ u\in R,~\alpha \in (-\frac{1}{2},\frac{1}{2}%
),~H=\alpha +1/2,  \label{1rep}
\end{equation}%
\textit{is the standard two-sided fBm }$W^{H}$\textit{.}

The proof and references within about other representations for $W^{H}$~can
be found in Section 2.

The particular case of (\ref{1rep}) with $\alpha \in (0,\frac{1}{2})$ and $%
a=b~$was noted in \cite{Pflug}, however, to the best of our knowledge the
case $\alpha \in (-\frac{1}{2},0)~$for the cusp signal has not been explored
so far and this paper fills the gap in the existing theory.

Further we discuss the model (\ref{swn}) with the signal of the form (\ref%
{signal}) and prove the asymptotic results$~$extending the aforementioned
results from \cite{Chernoyar} to the case $H=\alpha +\frac{1}{2}\in (0,1)$.

We denote the likelihood ratio (i.e. the Radon-Nykodim density, expressed in
terms of $X$, \cite{1}) as%
\begin{equation}
L(\theta ,X^{T}):=\exp \{\frac{1}{\varepsilon ^{2}}\int_{0}^{T}S(t,\theta
)dX_{t}-\frac{1}{2\varepsilon ^{2}}\int_{0}^{T}(S(t,\theta ))^{2}dt\}
\label{LikeSignal}
\end{equation}%
Our next main result is about the limit distributions for MLE\ $\hat{\theta}%
_{n}~$\ \eqref{LikeSignal} and general BE $\tilde{\theta}_{\varepsilon }^{.}$%
,\ under the assumption that the prior $q(\theta ),~\theta \in \Theta
=(\theta _{1},\theta _{2})~$\ is a continuous positive function, including
the Pitman-type estimate with noninformative prior 
\begin{equation*}
\tilde{\theta}_{\varepsilon }^{P}=\frac{\int_{_{_{\theta _{1}}}}^{\theta
_{2}}u\,L(u,X^{T})du}{\int_{_{_{\theta _{1}}}}^{\theta _{2}}\,L(u,X^{T})du}.
\end{equation*}%
For emphasising the dependence of the expected values and distributions on
an unknown parameter $\theta \ $below$\ $we will use the notation $\mathbf{E}%
_{\theta }\left( .\right) $~and $P_{\theta }^{.}(.).$

\textbf{Theorem 2.} \textit{Let (\ref{swn}) the cusp signal be defined as (%
\ref{signal}) with (\ref{cusp1}) such that }$\alpha \in (-1/2,1/2),$\textit{%
\ }$H=\alpha +1/2\in (0,1)$\textit{. Then as }$\varepsilon \rightarrow 0$%
\begin{equation}
\frac{(\hat{\theta}_{\varepsilon }-\theta )}{(\varepsilon /\Gamma _{\alpha
})^{\frac{1}{H}}}\overset{d}{\longrightarrow }\xi _{H},~\lim_{\varepsilon
\rightarrow 0}\mathbf{E}_{\theta }\left( \frac{\hat{\theta}_{\varepsilon
}-\theta }{(\varepsilon /\Gamma _{\alpha })^{\frac{1}{H}}}\right) ^{2}=%
\mathbf{E}\left( \xi _{H}^{2}\right) <\infty ,  \label{Th 2 st 2}
\end{equation}

\begin{equation}
\frac{(\tilde{\theta}_{\varepsilon }-\theta )}{(\varepsilon /\Gamma _{\alpha
})^{\frac{1}{H}}}\overset{d}{\longrightarrow }\zeta _{H},~\lim_{\varepsilon
\rightarrow 0}\mathbf{E}_{\theta }\left( \frac{\tilde{\theta}_{\varepsilon
}-\theta }{(\varepsilon /\Gamma _{\alpha })^{\frac{1}{H}}}\right) ^{2}=%
\mathbf{E}\left( \zeta _{H}^{2}\right) <\infty .  \label{Th 2 st1}
\end{equation}

\textit{Moreover, for any estimator }$\theta _{\varepsilon }=\theta
_{\varepsilon }(X^{T})$\textit{\ }%
\begin{equation}
\mathop{\underline{\lim}}\limits_{\delta \rightarrow 0}\mathop{\underline{%
\lim}}\limits_{\varepsilon \rightarrow 0}\sup_{\left\vert \theta -\theta
_{0}\right\vert \leq \delta }\mathbf{E}_{\theta }\left( \frac{\theta
_{\varepsilon }-\theta }{(\varepsilon /\Gamma _{\alpha })^{\frac{1}{H}}}%
\right) ^{2}\geq \mathbf{E}\left( \zeta _{H}^{2}\right) .  \label{Th 2 st3}
\end{equation}

The proof of Theorem 2 along with some discussions is presented in Section 3.

According to Theorem 2 both estimators $\hat{\theta}_{\varepsilon }$ and $%
\tilde{\theta}_{\varepsilon }$~have the same rate of convergence $%
\varepsilon ^{^{\frac{1}{H}}}\ $which is known to be the best possible rate.
It is natural to make a comparison of the properties of these estimators,
however, the analytical tools for studying of functionals of fBm are very
limited. In Section 4 we have included the series of simulation results to 
\textit{illustrate }properties of the limit random variables$~\zeta _{H}~$%
and $\xi _{H}~\,$for $H\in \lbrack 0.3,1);$ these results demonstrate that
the ratio of variances $\mathbf{E}\left( \xi _{H}^{2}\right) /\mathbf{E}%
\left( \zeta _{H}^{2}\right) ~$monotonically decreases from (approximately) $%
1.4.~$to $1$ as $H~\ $increases from $0.3$ to $1$.

\textbf{2. Representations of fBm }$W^{H},~H\in (0,1)$\textbf{.}

One can check (e.g. using Mathematica$^{@}$, version 11.0) that $\Gamma
_{\alpha }^{2}$ \ defined in (\ref{GGa}) can be expressed in terms of Gamma
function~$\Gamma \lbrack x]~$for any$~\alpha \in (-\frac{1}{2},0)\cup (\frac{%
1}{2},1)$ 
\begin{equation*}
\Gamma _{\alpha }^{2}=\frac{\sqrt{\pi }\Gamma \lbrack 1+\alpha ]}{2^{2\alpha
+1}\Gamma \lbrack 3/2+\alpha ]}(\sec [\pi \alpha ](a^{2}+b^{2})-2ab).
\end{equation*}%
For $\alpha \in (0,\frac{1}{2}),\ $see other equivalent representations in 
\cite{Prakasa}, p. 79 and \cite{1}, p. 306.

\textbf{Proof of Theorem 1. }The stochastic process $Y_{u}^{H}$ is an
Ito-type integral of a deterministic integrand and, obviously, $Y_{u}^{H}$ \
is a Gaussian process, \textbf{\ }$Y_{0}^{H}=0,\ EY_{u}^{H}=0.~$Thus for
verifying the statement of Theorem 1 we only need to find $%
E|Y_{u}^{H}-Y_{s}^{H}|^{2}$~for $u>s.$

Using the isometry property of Ito integrals and then the substitution $%
y=z+s~$\ we have%
\begin{eqnarray*}
E|Y_{u}^{H}-Y_{s}^{H}|^{2} &=&\Gamma _{\alpha }^{-2}\int_{-\infty }^{\infty
}(q_{\alpha }(y-u)-q_{\alpha }(y-s))^{2}dy \\
&=&\Gamma _{\alpha }^{-2}\int_{-\infty }^{\infty }(q_{\alpha
}(z-(u-s))-q_{\alpha }(z))^{2}dz.
\end{eqnarray*}%
Making the substitution $z=(u-s)x~$\ and using the identities 
\begin{equation*}
q_{\alpha }((u-s)x-(u-s))=|u-s|^{\alpha }q_{\alpha }(x-1),~q_{\alpha
}((u-s)x)=|u-s|^{\alpha }q_{\alpha }(x)
\end{equation*}%
we obtain 
\begin{equation*}
E|Y_{u}^{H}-Y_{s}^{H}|^{2}=\left\vert u-s\right\vert ^{2\alpha +1}\Gamma
_{\alpha }^{-2}\int_{-\infty }^{\infty }(q_{\alpha }(x-1)-q_{\alpha
}(x))^{2}dx=\left\vert u-s\right\vert ^{2H}
\end{equation*}%
where $H=\alpha +\frac{1}{2}.~$The proof is completed.

\textbf{Remark 1. }

For the symmetric case $a=b$ representation (\ref{1rep}) was obtained in 
\cite{Pflug}.

If $~a=1,~b=0~$it is equivalent to the Mandelbrot-Van Ness representation:%
\begin{equation*}
Y_{u}^{H}=\Gamma _{\alpha }^{-1}(\int_{0}^{u}(u-y)^{\alpha }dW\left(
y\right) +\int_{-\infty }^{0}((u-y)^{\alpha }-(-y)^{\alpha })dW\left(
y\right) ),u\in R,
\end{equation*}%
see \cite{M-N}.

The Muravlev's representation \cite{Mur} in terms of Ornstein-Uhlenbeck
processes is a consequence of the Mandelbrot-Van Ness representation;
actually, the latter can be considered as a consequence of Kolmogorov's
representation \cite{Kol} for fBm in terms of a Fourier transform of a
Gaussian random field. There exist other representations for the fBm, not
directly connected to (\ref{1rep}), see e.g. Norros et al. \cite{Nor}.

\textbf{Remark 2. } Theorem 6.2.1 from \cite{1}, applied to a particular
case of cusp densities (\ref{cusp1}),\ contains a representation for the
limit of normalized likelihood ratio process (NLRP) in the form of
stochastic integrals of cusp-type functions. Interestingly, the connection
of such stochastic processes to the fBm had not been discussed in \cite{1}
at all. Now using (\ref{1rep}) one can easily check that in \cite{1} the
Gaussian component in the limit of NLRP for the case under consideration is
nothing else but the fBm $W^{H},~H\in (\frac{1}{2},1)$. \footnote{%
This fact has not been clarified in the literature before.}

\textbf{3. Cusp-type signals in "signal plus white noise" model}

Further we use the following notation for the NLRP 
\begin{equation*}
Z_{T}\left( u,\varepsilon \right) :=\frac{L(\theta -\varphi _{\varepsilon
}u,X^{T})}{L(\theta ,X^{T})}
\end{equation*}%
\begin{equation*}
=\exp \{\varepsilon ^{-2}\int_{0}^{T}[S(t,\theta -\varphi _{\varepsilon
}u)-S(t,\theta )]dX_{t}-\frac{\varepsilon ^{-2}}{2}\int_{0}^{T}(S^{2}(t,%
\theta -\varphi _{\varepsilon }u)-S^{2}(t,\theta ))dt\},
\end{equation*}%
where we assume $u\in U_{\varepsilon }:=(\frac{\theta _{1}-\theta }{\varphi
_{\varepsilon }},\frac{\theta _{2}-\theta }{\varphi _{\varepsilon }})~$and
set%
\begin{equation*}
\varphi _{\varepsilon }:=(\varepsilon /\Gamma _{\alpha })^{\frac{1}{H}},
\end{equation*}%
thus $\varepsilon ^{-2}\varphi _{\varepsilon }^{2H}\Gamma _{\alpha }^{2}=1$.
Having chosen $\varphi _{\varepsilon },$ this way we obtain the
representations for the limit distributions of $\hat{\theta}_{\varepsilon }$%
\ and $\tilde{\theta}_{\varepsilon }~$identical to these in (\ref{limit1})
and (\ref{MLE limit}). The limit distributions in \cite{Chernoyar} can be
transformed to (\ref{Th 2 st 2}) and (\ref{Th 2 st1}) after properly
adjusting the normalising factor $\varphi _{\varepsilon }$.

The proof of Theorem 2 is based on the properties of NLRP $Z_{T}\left(
u,\varepsilon \right) $\textit{.~ }First, we prove the convergence$~$of
marginal distributions of\textit{\ }$Z_{T}\left( u,\varepsilon \right) $ to $%
Z_{u}^{(H)}=\exp \{W_{u}^{H}-\frac{|u|^{2H}}{2}\}$ as $\varepsilon
\rightarrow 0,\ u\in R,~$ $H=\alpha +1/2.$

\textbf{Proposition 1.}\textit{\ Assume (\ref{swn}) holds where}%
\begin{equation*}
S(t,\theta )=q_{\alpha }(t-\theta )+h(t,\theta ),~\alpha \in (-\frac{1}{2},%
\frac{1}{2}),
\end{equation*}%
$h(t,\theta )\in C^{1,1},\,0<\theta _{1}<\theta <\theta _{2}<T.$

\textit{Then the marginal distributions of }$\{Z_{T}\left( u,\varepsilon
\right) ,\ u\in U_{\varepsilon }\}$\textit{\ converge to marginal
distributions of }$\{Z_{u}^{(H)},\ u\in R\},~H=\alpha +\frac{1}{2}.$

\textbf{Proof. }Using the equation~$dX_{t}=S(t,\theta )dt+\varepsilon dw_{t}$
we have 
\begin{eqnarray*}
Y_{\varepsilon }(u) &:&=\log (Z_{T}\left( u,\varepsilon \right) ) \\
&=&\varepsilon ^{-1}\int_{0}^{T}[S(v,\theta -\varphi _{\varepsilon
}u)-S(v,\theta )]dw_{v}-\frac{\varepsilon ^{-2}}{2}\int_{0}^{T}[S(v,\theta
-\varphi _{\varepsilon }u)-S(v,\theta )]^{2}dv \\
&:&=A_{\varepsilon }(u)-B_{\varepsilon }(u).
\end{eqnarray*}%
First we show, as $\varepsilon \rightarrow 0,~$ the deterministic part of
this decomposition 
\begin{equation*}
B_{\varepsilon }(u)\rightarrow \frac{|u|^{2H}}{2},\ u\in R,
\end{equation*}%
and the stochastic integral part%
\begin{equation*}
E_{\theta }|A_{\varepsilon }(u)-A_{\varepsilon }(s)|^{2}\rightarrow
E|W_{u}^{H}-W_{s}^{H}|^{2}=|u-s|^{2H},\,\ u\in R,~s\in R.
\end{equation*}%
The last convergence is equivalent to the convergence of covariance
functions of $A_{\varepsilon }(u)~$to $W^{H}~$\ and since\ $A_{\varepsilon
}(u)$ is a Gaussian process this will imply the result.

We have 
\begin{equation*}
B_{\varepsilon }(u)=\varepsilon ^{-2}\int_{0}^{T}[q_{\alpha }(v-\theta
-\varphi _{\varepsilon }u)-q_{\alpha }(v-\theta )+\delta _{\varepsilon
}(v,u)]^{2}dv,
\end{equation*}%
where 
\begin{equation*}
\delta _{\varepsilon }(v,u):=h(v,\theta -\varphi _{\varepsilon
}u)-h(v,\theta )=\varphi _{\varepsilon }D(v,\theta ,u)(1+o(1))=o(\varphi
_{\varepsilon }).
\end{equation*}%
Since $h(v,\theta )\in C^{1,1}$ one can easily see that the total input of $%
\delta _{\varepsilon }(t,u)$ to $B_{u}(\varepsilon )$ is of order $o(1)~$%
because%
\begin{equation*}
\varepsilon ^{-2}\int_{0}^{T}(\delta _{\varepsilon }(v,u))^{2}dv\leq
\varepsilon ^{-2}\varphi _{\varepsilon }^{2}\max_{v}|D(v,\theta ,u)|^{2}T
\end{equation*}%
and $\varepsilon ^{-2}\varphi _{\varepsilon }^{2}=\varepsilon ^{-2}\varphi
_{\varepsilon }^{2H}\varphi _{\varepsilon }^{2-2H}=o(1)$.~Hence, we obtain
with the substitution $v=\theta +\varphi _{\varepsilon }t~~$ 
\begin{equation*}
B_{\varepsilon }(u)=\varepsilon ^{-2}\int_{0}^{T}[q_{\alpha }(v-\theta
-\varphi _{\varepsilon }u)-q_{\alpha }(v-\theta )]^{2}dv+o(1)
\end{equation*}%
\begin{equation*}
=\varepsilon ^{-2}\varphi _{\varepsilon }^{2\alpha +1}\int_{-\theta /\varphi
_{\varepsilon }}^{(T-\theta )/\varphi _{\varepsilon }}[q_{\alpha
}(t-u)-q_{\alpha }(t)]^{2}dt+o(1)
\end{equation*}%
(recall $\varepsilon ^{-2}\varphi _{\varepsilon }^{2\alpha +1}\Gamma
_{\alpha }^{2}=1,~T>\theta _{2}\geq \theta \geq \theta _{1}>0$)%
\begin{equation*}
\rightarrow \Gamma _{\alpha }^{-2}\int_{-\infty }^{\infty }[q_{\alpha
}(t-u)-q_{\alpha }(t)]^{2}dt=|u|^{2H},\ u\in R.
\end{equation*}%
Due to the isometry of stochastic integrals we obtain 
\begin{eqnarray*}
&&E_{\theta }|A_{\varepsilon }(u)-A_{\varepsilon }(s)|^{2} \\
&=&\varepsilon ^{-2}\int_{0}^{T}[q_{\alpha }(v-\theta -\varphi _{\varepsilon
}u)-q_{\alpha }(v-\theta -\varphi _{\varepsilon }s)+\delta _{\varepsilon
}(t,u)-\delta _{\varepsilon }(t,s)]^{2}dv
\end{eqnarray*}%
\begin{equation*}
=\varepsilon ^{-2}\varphi _{\varepsilon }^{2\alpha +1}\int_{-\theta /\varphi
_{\varepsilon }}^{(T-\theta )/\varphi _{\varepsilon }}[q_{\alpha
}(t-u)-q_{\alpha }(t-s)]^{2}dt+o(1)\rightarrow 
\end{equation*}%
\begin{equation*}
\Gamma _{\alpha }^{-2}|u-s|^{2H}\int_{-\infty }^{\infty }[q_{\alpha
}(t-u)-q_{\alpha }(t)]^{2}dt=|u-s|^{2H}
\end{equation*}%
as $\varepsilon \rightarrow 0~$. This completes the proof.%
\begin{equation*}
\end{equation*}%
\textbf{Remark 3.} Proposition 1 represents the extension of Lemma 1 from 
\cite{Chernoyar} to the case of asymmetric cusp signals and the full range $%
\alpha \in (-\frac{1}{2},\frac{1}{2})~$albeit with the essentially shortened
proof. The extension~is achieved due to the exploitation of the fact that
the convergence of marginal distributions of any Gaussian process $%
Y_{\varepsilon }(u)~$is equivalent to the convergence of its covariance
functions $Cov(Y_{\varepsilon }(u),Y_{\varepsilon }(s))$ as ~\ $\varepsilon
\rightarrow 0.$

\textbf{Proposition 2. }\textit{Under conditions of Proposition 1 }$\ $the
process $\{Z_{T}\left( u,\varepsilon \right) ,\ u\in \lbrack a,b]\}~$\textit{%
converges weakly to }$\{Z_{u}^{(H)},~u\in \lbrack a,b]\},~H=\alpha +1/2\in
(0,1)~$\textit{in the space of continuous functions }$C[a,b]$ \textit{for
any finite interval }$[a,b].$

Note that our proof of Theorem 2 will consist in verifying conditions of the
fundamental Theorems 1.10.1 and 1.10.2 \ from \cite{1} and it will not rely
on Proposition 2. However, we believe that it is useful to make a simple
demonstration how the fBm appears as a process with trajectories in $C[a,b]$
in the limit. At the same time we would like to stress that the enhancement
from $C[a,b]$ to $C(-\infty ,\infty )$ requires some extra conditions,~see
Remark 4 and condition (\ref{ccc1}) below.

\textbf{Proof. }We will show the convergence of the\textit{\ continuous
Gaussian} process $Y_{\varepsilon }(u)=\log Z_{T}(u,\varepsilon )~$to $%
Y(u):=\log Z_{u}^{(H)}.~$This requires verifying the technical condition for
convergence in $C[a,b],$ see e.g. Theorem 4.3 and 4.4. in \cite{Prakasa},
also see Gikhman and Skorokhod \cite{GihSk}. According to these references
it is sufficient\ to check that there exist the constants $\ p>0$ and $q>1~$%
such that all $u~$and $s~$\ from any interval $\,[a,b]\in R$%
\begin{equation*}
E_{\theta }|Y_{\varepsilon }(u)-Y_{\varepsilon }(s)|^{p}\leq C|u-s|^{q},
\end{equation*}%
where $C\,\ $is a generic constant, which does not depend on $\varepsilon
,~u~$\ and $s$.~Using the notation for $A_{\varepsilon }(u)$ and $%
B_{\varepsilon }(u)$ introduced in the proof of Proposition 1 above we
obtain 
\begin{equation}
E_{\theta }|Y_{\varepsilon }(u)-Y_{\varepsilon }(s)|^{p}\leq C(E_{\theta
}|A_{\varepsilon }(u)-A_{\varepsilon }(s)|^{p}+|B_{\varepsilon
}(u)-B_{\varepsilon }(s)|^{p})  \label{w1}
\end{equation}%
Since $A_{\varepsilon }(u)$ is a Gaussian process we have 
\begin{equation}
E_{\theta }|A_{\varepsilon }(u)-A_{\varepsilon }(s)|^{p}=C(E|A_{\varepsilon
}(u)-A_{\varepsilon }(s)|^{2})^{p/2}.  \label{w2}
\end{equation}%
Using the arguments from the proof of Proposition 1 above and the inequality 
$(x+y)^{2}\leq 2(x^{2}+y^{2})$ we obtain%
\begin{equation*}
E_{\theta }|A_{\varepsilon }(u)-A_{\varepsilon }(s)|^{2}=\varepsilon
^{-2}\int_{0}^{T}[q_{\alpha }(v-\theta +\varphi _{\varepsilon }u)-q_{\alpha
}(v-\theta +\varphi _{\varepsilon }s)+\delta _{\varepsilon }(t,u)-\delta
_{\varepsilon }(t,s)]^{2}dv
\end{equation*}%
\begin{equation*}
\leq 2\varepsilon ^{-2}\varphi _{\varepsilon }^{2H}\int_{-\theta /\varphi
_{\varepsilon }}^{(T-\theta )/\varphi _{\varepsilon }}[q_{\alpha
}(t-u)-q_{\alpha }(t-s)]^{2}dt+2\varepsilon ^{-2}\varphi _{\varepsilon
}^{2}T\max_{v}|D(v,\theta ,u)|^{2}|u-s|^{2}
\end{equation*}%
(substituting $y=(u-s)x$) 
\begin{equation*}
\leq 2\Gamma _{\alpha }^{-2}|t-s|^{2\alpha +1}|\int_{[-s\ -\ \theta /\varphi
_{\varepsilon }]/(u-s)}^{[-s\ +\ (T-\theta )/\varphi _{\varepsilon
}]/(u-s)}[q_{\alpha }(v-1)-q_{\alpha }(v)]^{2}dt|+C|u-s|^{2}.
\end{equation*}%
Hence,%
\begin{equation}
E_{\theta }|A_{\varepsilon }(u)-A_{\varepsilon }(s)|^{2}\leq
C(|u-s|^{2\alpha +1}+|u-s|^{2})  \label{w3}
\end{equation}%
Since $|B_{\varepsilon }(u)-B_{\varepsilon }(s)|^{p}\leq C|u-s|^{p},~$now (%
\ref{w1}),~(\ref{w2}) and (\ref{w3})~imply%
\begin{equation*}
E_{\theta }|Y_{\varepsilon }(u)-Y_{\varepsilon }(s)|^{p}\leq
C(|u-s|^{(2\alpha +1)p/2}+|u-s|^{p}).
\end{equation*}%
Choosing $p$ large enough such that $q=(2\alpha +1)p/2>1$\ then for $u\in
\lbrack a,b]$ and $s\in \lbrack a,b]$ we obtain 
\begin{equation*}
E_{\theta }|Y_{\varepsilon }(u)-Y_{\varepsilon }(s)|^{p}\leq
C|u-s|^{(2\alpha +1)p/2}(1+|u-s|^{p(\frac{1}{2}-\alpha )})\leq C|u-s|^{q}.
\end{equation*}%
This completes the proof.%
\begin{equation*}
\end{equation*}%
\textbf{Remark 4. }The extension from $C[a,b]$ to $C(-\infty ,\infty )$ in
Proposition 2 can not be done without verifying some extra conditions, see
the counterexample in \cite{YK 1998}, Remark 4.2, p. 161. and condition (\ref%
{ccc1}) below.

\textbf{Proof of Theorem 2.} The detailed exposition of the technique
required for the proof can be found in \cite{1}, \cite{Kut book} or in \cite%
{Chernoyar}. Hence, in addition to Proposition 1, following a$\ $%
well-trodden path we need to make the following steps.

\textbf{Step 1. }As clarified in \cite{Chernoyar} to apply Theorem 1.10.1
and 1.10.2 from \cite{1} to the continuous time model (\ref{swn}) we need
first to prove of convergence of marginal distributions (which is done above
in Proposition 1) and also show that there exists $C>0$ such that

\begin{equation}
\mathbf{E}_{\theta }\left( \sqrt{Z_{T}\left( u,\varepsilon \right) }\right)
\leq \exp \{-C|u|^{2H}\}.  \label{ccc1}
\end{equation}

Under our assumptions for the general cusp (\ref{cusp1}), the inequality (%
\ref{ccc1}) can be proved by mimicking the proofs of Lemma 2 and 3 from \cite%
{Chernoyar}. In particular, at first we show that there exists $C>0~$such
that for $u\in U_{\varepsilon }$%
\begin{equation*}
\int_{0}^{T}(S^{2}(v,\theta -\varphi _{\varepsilon }u)-S^{2}(v,\theta
))dv\geq C|\varphi _{\varepsilon }u|^{2H}
\end{equation*}%
and then noting that 
\begin{equation*}
\mathbf{E}_{\theta }\left( \sqrt{Z_{T}\left( u,\varepsilon \right) }\right)
=\exp \{-\frac{1}{8\varepsilon ^{2}}\int_{0}^{T}(S^{2}(v,\theta -\varphi
_{\varepsilon }u)-S^{2}(v,\theta ))dv\}
\end{equation*}%
conclude that (\ref{ccc1}) holds.

\textbf{Step 2. }Accordingly to Theorem 1.10.1 and 1.10.2 from \cite{1} we
need to show that for $m>0~$and any $u\in U_{\varepsilon },~s\in
U_{\varepsilon }~$there exists $\beta >1$ such that 
\begin{equation}
\mathbf{E}_{\theta }|(Z_{T}\left( u,\varepsilon \right)
)^{1/(2m)}-(Z_{T}\left( s,\varepsilon \right) )^{1/(2m)}|^{2m}\leq
C|u-s|^{\beta },  \label{cc2}
\end{equation}%
where $C$ is a generic positive constant.~This can be done by showing 
\begin{eqnarray*}
&&\mathbf{E}_{\theta }|(Z_{T}\left( u,\varepsilon \right)
)^{1/(2m)}-(Z_{T}\left( s,\varepsilon \right) )^{1/(2m)}|^{2m}\leq \\
&\leq &C(\varepsilon ^{-2}\int_{0}^{T}[S(v,\theta -\varphi _{\varepsilon
}u)-S(v,\theta -\varphi _{\varepsilon }s)]^{2}dv)^{m}
\end{eqnarray*}%
in the lines of the proof of Proposition 2 above. Then using the estimates
obtained in the proofs of Propositions 1 and 2, we can easily check that 
\begin{equation*}
(\varepsilon ^{-2}\int_{0}^{T}[S(v,\theta -\varphi _{\varepsilon
}u)-S(v,\theta -\varphi _{\varepsilon }s)]^{2}dv)^{m}\leq C|u-s|^{2mH}
\end{equation*}%
then~choosing $m$ such that $2mH>1$,$\ $we obtain (\ref{cc2}).

Finally, based on \ (\ref{ccc1})\ and \ (\ref{cc2}) the tightness of the
family of the distributions of $\ Z_{T,\ \theta }\left( \varepsilon \right)
\ $process$\ $can be proved in the following sense,\ see \cite{1}, Chapter
1, \ Theorem 5.1 and Remark 5.1. For any $N>0\ $\ and any compact set $K\ \ $%
there exist $M_{0}\ $and $c_{N}\ $such that for any $M>\ M_{0}\ \ $and all $%
\varepsilon <\varepsilon (N,\ K)$,%
\begin{equation*}
\underset{\theta \in K}{\sup }P_{\theta }^{.}\left( \underset{|u|>M}{\sup }%
Z_{T}\left( u,\varepsilon \right) >M^{-N}\right) \leq c_{N}M^{-N}.
\end{equation*}

This completes the proof.

\textbf{Remark 5. }The uniqueness with probability one of the random
variable $\xi _{H}$ can be shown also using the standard arguments related
to the continuous mapping theorem for argmax functionals, see \cite%
{KimPollard}.

\textbf{4. Simulations results for the densities and variances of the limit
distributions.}

To apply results of Theorem 2 for constructing asymptotic confidence
intervals for $\theta ~$it is desirable to know densities of $\zeta _{H}~$%
and $\xi _{H}~$but to our best knowledge there are no general analytical or
numerical methods for this purpose. The difficulty is due to the fact that
for $H\neq 1~$and $H\neq \frac{1}{2}~$ the fBm $W_{u}^{H}~$ is neither a
Markov process nor a semimartingale, rendering the standard tools of Markov
theory and stochastic analysis are not applicable, at least directly. Some
general properties of $\zeta _{H},~H\in (0,1),$ has been obtained in \cite%
{NK}, \cite{NKL} with the help of the measure transformation technique.

It is well known that at the boundary point $H=1$ both $\zeta _{1}~$and $\xi
_{1}$ have a standard normal distribution and so $Var(\zeta _{1})=Var(\xi
_{1})=1.$ Besides$~$\ the case $H=1$ there is only one explicit analytical
result for the density of $\xi _{\frac{1}{2}}\ $obtained in \cite{Yao}, \cite%
{Shepp}: 
\begin{equation}
P(|\xi _{\frac{1}{2}}|>t)=(t+5)\Phi (-\frac{\sqrt{t}}{2})-\sqrt{\frac{2t}{%
\pi }}e^{-\frac{t}{8}}-3e^{t}\Phi (-\frac{3\sqrt{t}}{2}),  \label{yao}
\end{equation}%
where $\Phi (t)$ is a standard normal distribution. This result implies 
\begin{equation}
Var(\xi _{\frac{1}{2}})=26,  \label{Tere26}
\end{equation}%
the latter firstly was obtained in \cite{Tere}, see also \cite{1}.

The analytical form for the density of $\zeta _{\frac{1}{2}}$ is still
unknown but in \cite{RS} (see also \cite{NKL}, \cite{NK}) it was shown 
\begin{equation}
Var(\zeta _{\frac{1}{2}})=16Zeta[3]\thickapprox 19.23,
\label{explicit dzeta}
\end{equation}%
where $Zeta[k]$ is the Euler-Riemann's zeta-function.

To simulate $\zeta _{H}$ and $\xi _{H}~$for arbitrary $H~$we truncated the
integration range $u\in (-\infty ,\infty )$ to $u\in \lbrack -T,T]$ and then
simulated discretised fBm trajectories $\left\{ W_{u_{j}}^{H,(i)}\right\}
_{j=-m}^{m},~u_{j}\in \left\{ jT/m\right\} _{j=-m}^{m},j\in \mathbb{Z},$
based on the Wood-Chan's algorithm \cite{AnWood}. Note that errors due
discretisation of fBm trajectories are of order $O(m^{-H})~$\ and so they
could be significant when $H$ is small even with relatively large $m=2^{19}~$%
(see in \cite{B} some results about the rate of convergence of
max-functionals~of $\left\{ W_{u_{j}}^{H,(i)}\right\} _{j=-m}^{m}$ to the
limit). That is the reason why we decided not to include simulation results
for values $H<0.3~$where we did observe significant errors. Potentially,
more accurate results can be obtained with values $m=2^{20}$ or higher but
this would take much more computational time which was not affordable even
for high performance computers available to us.

For $i$-th simulation, $i=1,\ldots ,N$, we approximate $\zeta _{H}$ by 
\begin{equation}
\widehat{\zeta }_{H}^{\;(i)}=\frac{\sum_{j=-m}^{m}u_{j}\;\exp \left\{
W_{u_{j}}^{H,(i)}-\frac{1}{2}|u_{j}|^{2H}\right\} \;w(u_{j})}{%
\sum_{j=-m}^{m}\exp \left\{ W_{u_{j}}^{H,(i)}-\frac{1}{2}|u_{j}|^{2H}\right%
\} \;w(u_{j})}~,  \label{BE:limit:sim}
\end{equation}%
where $w(u_{j})$ are trapezoidal rule weights, and approximate $\xi _{H}$ by 
\begin{equation}
\widehat{\xi }_{H}^{\;(i)}=\argmax_{\scriptscriptstyle u_{j}\in \left\{
jT/m\right\} _{j=-m}^{m},j\in \mathbb{Z}}\left\{ W_{u_{j}}^{H,(i)}-\frac{1}{2%
}|u_{j}|^{2H}\right\} \;.  \label{MLE:limit:sim}
\end{equation}
Sample variances of limit distributions $\zeta _{H}$ and $\xi _{H}$, denoted
by $\widehat{{Var}}[\zeta _{H}]$ and $\widehat{{Var}}[\xi _{H}]$
respectively, are depicted in Table \ref{zeta:xi:var} and Figure \ref%
{Fig:varxi_varzeta}. The sample variance reported for $H\geq 0.3$ are
calculated based on the random variables $\left\{ \zeta _{H}^{\;(i)}\right\}
_{i=1}^{N}$ and $\left\{ \xi _{H}^{\;(i)}\right\} _{i=1}^{N}$ simulated
using the setting $N=10^{7}$, $m=2^{19}$, $T=10^{5}$, $u\in \lbrack -T,T]$.
In \cite{NKL}, sample variance of $\left\{ \zeta _{H}^{\;(i)}\right\}
_{i=1}^{N}$ and $\left\{ \xi _{H}^{\;(i)}\right\} _{i=1}^{N}$ used to
simulate the limit distributions $\zeta _{H}$ and $\xi _{H}$ estimated using %
\eqref{BE:limit:sim} and \eqref{MLE:limit:sim} were reported for $H\in
\lbrack 0.4,0.91]$ using the setting $N=10^{6}$, $m=2^{18}$, $T=10^{5}$, $%
u\in \lbrack -T,T]$. In the current work, we report results for sample
variances of simulated random variables $\left\{ \zeta _{H}^{\;(i)}\right\}
_{i=1}^{N}$ and $\left\{ \xi _{H}^{\;(i)}\right\} _{i=1}^{N}$ for a wider
range of $H$, i.e., $H\in \lbrack 0.3,0.99]$, together with 10 times more
simulated trajectories, and two times more discrete points on either side of
zero while using the same truncation limit $[-10^{5},10^{5}]$.

\begin{figure}[!htb]
\centering
\par
\par
\includegraphics[width=12cm, height=8cm]
{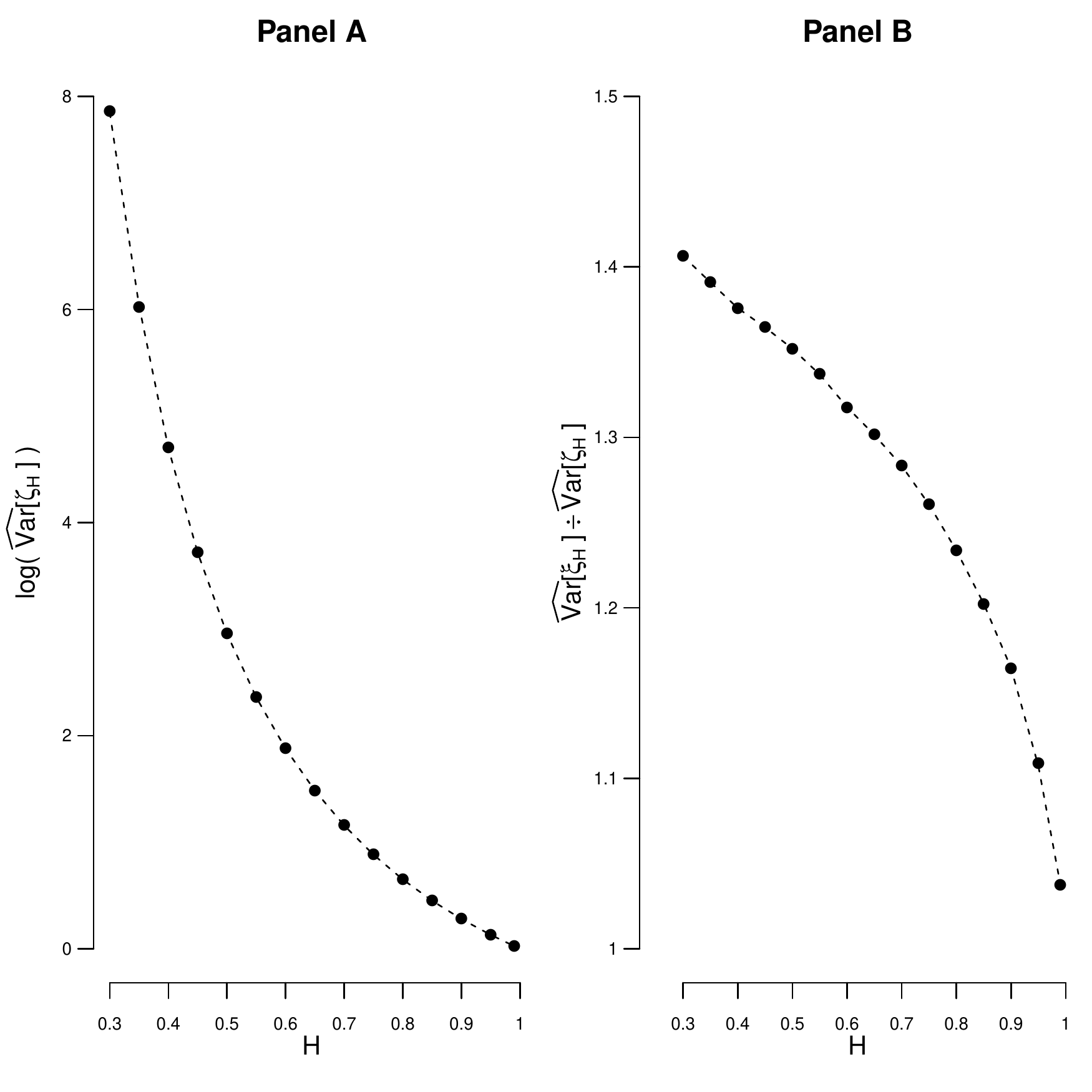}
\caption{\textbf{Panel A:} $\log(\widehat{{Var}}[\protect\zeta_{H}])$
against $H$. \textbf{Panel B:} $\widehat{{Var}}[\protect\xi_{H}] / \widehat{{%
Var}}[\protect\zeta_{H}]$ against $H$.}
\label{Fig:varxi_varzeta}
\end{figure}

For the case $H=0.5$~we obtained $\widehat{{Var}}[\zeta _{1/2}]=19.206$,
which is close to $16Zeta[3]\approx 19.23$. Table \ref{zeta:xi:var} also
depict the sample variance $\widehat{{Var}}[\xi _{H}] $ against the Hurst
parameter $H$. As expected, the $\widehat{{Var}}[\zeta _{H}]$ is smaller
than $\widehat{{Var}}[\xi _{H}]$. Figure \ref{Fig:varxi_varzeta} illustrates
this point from a graphical perspective.

\begin{figure}[!htb]
\centering
\par
\par
\includegraphics[width=15cm,
height=12cm]{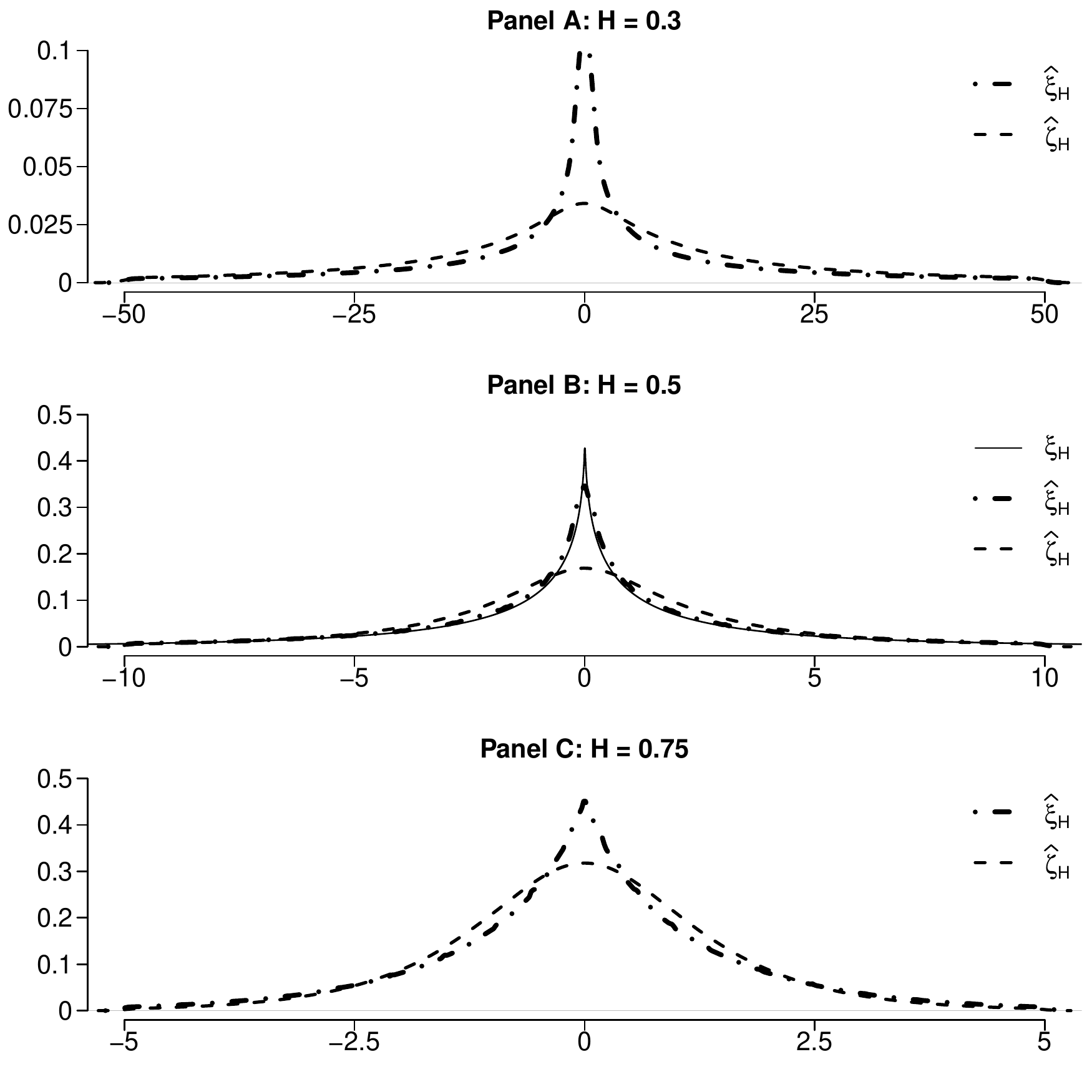}
\caption{$\protect\xi_H$: True probability function of $\protect\xi_H$ for $%
H=0.5$. $\widehat{\protect\xi}_{H}$: Kernel density approximation of
probability density function from $\widehat{\protect\xi}^{\;(i)}_{H}$. $%
\widehat{\protect\zeta}_{H}$: Kernel density approximation of probability
density function from $\widehat{\protect\zeta}^{\;(i)}_{H}$.}
\label{Fig:Density:Xi:Zeta}
\end{figure}

Figure \ref{Fig:Density:Xi:Zeta} depicts the approximate probability density
function of $\zeta _{H}$ and $\xi _{H}$ obtained by applying kernel density
smoothing on the simulated random variables $\left\{ \zeta
_{H}^{\;(i)}\right\} _{i=1}^{N}$ and $\left\{ \xi _{H}^{\;(i)}\right\}
_{i=1}^{N}$. 
\begin{table}[tbph]
\caption{Tabulation of sample variance, $\widehat{{Var}}[\protect\zeta _{H}]$
and $\widehat{{Var}}[\protect\xi _{H}]$, estimated empirically from $\left\{ 
\protect\zeta _{H}^{\;(i)}\right\} _{i=1}^{N}$ and $\left\{ \protect\xi %
_{H}^{\;(i)}\right\} _{i=1}^{N}$ that are calculated from $N=10^{7}$
simulated fBm trajectories, each simulated at $n=2^{19}$ equally spaced
discretization points on either side of zero spanning the interval $%
[-T,T],T=10^{5}$ for various values of Hurst's parameter $H$.}
\label{zeta:xi:var}
\newpage
\par
\begin{center}
\scalebox{0.75}{
\begin{tabular}{cllllllll}
  \hline 
\\
$H$
&\multicolumn{1}{c}{0.30}
&\multicolumn{1}{c}{0.35}      
&\multicolumn{1}{c}{0.40}      
&\multicolumn{1}{c}{0.45}     
&\multicolumn{1}{c}{0.50}     
&\multicolumn{1}{c}{0.55}
&\multicolumn{1}{c}{0.60} \\
 
 \hline   \hline 
\\
$\widehat{\Var}[\zeta_{H}]$ 
&2587.91
&411.22
&110.12
&41.17
&19.21
&10.58
&6.54 \\

$\widehat{\Var}[\xi_{H}]$
&3639.31
&572.05
&151.48 
&56.18
&25.97
&14.15 
&8.61  \\

  \hline 
\\
$H$
&\multicolumn{1}{c}{0.65}    
&\multicolumn{1}{c}{0.70}    
&\multicolumn{1}{c}{0.75}    
&\multicolumn{1}{c}{0.80}    
&\multicolumn{1}{c}{0.85}
&\multicolumn{1}{c}{0.90}    
&\multicolumn{1}{c}{0.95}    
&\multicolumn{1}{c}{0.99}    \\   

 \hline   \hline 
\\

$\widehat{\Var}[\zeta_{H}]$
&4.41
&3.18
&2.42
&1.91
&1.57
&1.32 
&1.14
&1.03\\

$\widehat{\Var}[\xi_{H}]$
&5.74 
&4.08
&3.04
&2.36
&1.88
&1.54
&1.26
&1.06\\

   \hline
\end{tabular}
}
\end{center}
\end{table}

\clearpage\newpage

\textbf{5. Conclusions.} \vskip 0.25cm

This paper presents some extensions of the results from \cite{Chernoyar}
where an estimation of a singularity point of a cusp-type signal in the
"signal plus white noise" was discussed. We demonstrated that when the
intensity of white noise $\varepsilon \rightarrow 0$ the limits of BE $%
\tilde{\theta}_{\varepsilon }^{.}$ and MLE $\hat{\theta}_{\varepsilon }$ are
expressed in terms of fBm $W^{H}$ for the full range of cusp-type signals
with $\alpha \in (-\frac{1}{2},\frac{1}{2})~$and correspondingly with the
Hurst parameter $H=\alpha +\frac{1}{2}\in (0,1).~$The simulations results
suggest that as $\varepsilon \rightarrow 0~$the $\ $limit of $Var(\hat{\theta%
}_{\varepsilon })/Var(\tilde{\theta}_{\varepsilon }^{.})$ $~$is a decreasing
function in the range $H\in \lbrack 0.3,1]~$showing about $40\%$ gain $%
\tilde{\theta}_{\varepsilon }$ over $\hat{\theta}_{\varepsilon }~~$\ for $%
H=0.3.$~and confirming the known result for $H=\frac{1}{2}~$where the
corresponding gain is about $35\%.~$

\textbf{6. Acknowledgements. }This paper was completed during the research
stay of the third author at Tokyo Metropolitan University (TMU). We thank
our colleagues from TMU for research assistance.

The authors also thank participants of Hokkaido Winter School on Finance and
Mathematics (February 2017) for discussions which helped to improve the
exposition of the paper.

The simulation results were obtained using the NCI National Facility in
Canberra, Australia, which is supported by the Australian Commonwealth
Government and on the eResearch HPCC at the University of Technology Sydney.
We gratefully acknowledge Dr. Timothy Gregory Ling for kindly providing us
with the \texttt{C++} code for simulation of fBm, and estimation of $\zeta
_{H}$ and $\xi _{H}$.

\end{document}